\documentclass[letterpaper]{article}

\usepackage{amsfonts}
\usepackage{amsmath}
\usepackage{amssymb}
\usepackage{latexsym}
\usepackage[letterpaper,pdftex=true,colorlinks,linkcolor=blue,citecolor=red,pdfstartview=FitH]{hyperref}
\providecommand{\U}[1]{\protect\rule{.1in}{.1in}}

\providecommand{\U}[1]{\protect\rule{.1in}{.1in}}

\newtheorem{teo}{Theorem}
\newtheorem{lem}{Lemma}

\newtheorem{com}{Comment}

\newcommand{\proof}{\textbf{Proof. }}

\begin{document}

\title{Riemann Hypothesis for Goss $t$-adic Zeta Function}
\date{}
\author{Javier Diaz-Vargas and Enrique Polanco-Chi.\\Universidad Aut\'{o}noma de Yucat\'{a}n, Facultad de Matem\'{a}ticas.\\Perif\'{e}rico Norte Tablaje 13615, 97119, M\'{e}rida, Yucat\'{a}n, M\'{e}xico.\\jdvargas@uady.mx, polanco89@live.com.mx.}
\maketitle

\begin{abstract}

In this short note, we give a proof of the Riemann hypothesis for
Goss $v$-adic zeta function $\zeta_{v}(s)$, when $v$ is a prime of
$\mathbb{F}_{q}[t]$ of degree one.

\end{abstract}

The Riemann Hypothesis says that the non-trivial complex zeros of the Riemann
zeta function all lie on a line Re$(s)=1/2$ in the complex plane. An analog
\cite{rhcpzf, bsffa, ffa} of this
statement was proved for the Goss zeta function for $\mathbb{F}_q[t]$, for
$q$ a prime by Wan \cite{rhcpzf} (see also \cite{rh}) and for general $q$
by Sheats \cite{rhgzf}. The proofs use the calculation of the Newton polygons associated
to the power series these zeta functions represent, and  their  slopes are calculated or
estimated by the degrees of power sums which make up the terms. In \cite{rh, rhgzf},
the exact degrees derivable (as noticed by Thakur) from incomplete work of Carlitz
were justified and the Riemann hypothesis was derived.

In this paper, we look at the $v$-adic Goss zeta function and prove
analog of the Riemann hypothesis, in case where $v$ is a prime of
degree one in $\mathbb{F}_q[t]$, using the valuation formulas for
the corresponding power sums given essentially (see below) in \cite{vpsgvzf}
and following the method of \cite{rh},
\cite[Sec. 5.8]{ffa}. We note that for $q=2$, this was already shown
by Wan \cite{rhcpzf} using the earlier calculation at the infinite
place. See below for the details. We also note that the situation
for higher degree $v$ is different \cite[Sec. 10]{vpsgvzf} and is not
fully understood.

We now give the relevant definitions and describe the results precisely.

In  the function field - number field analogy, we have
$A=\mathbb{F}_q[t]$, $K=\mathbb{F}_q(t)$,
$K_{\infty}=\mathbb{F}_q((1/t))$ and $\mathbb{C}_{\infty}$, the
completion of algebraic closure of $K_{\infty}$ as analogs of
$\mathbb{Z}$, $\mathbb{Q}$, $\mathbb{R}$ and $\mathbb{C}$
respectively. We consider $v$-adic situation where  $v$ is an
irreducible polynomial of $A$, and $K_v$ the completion  of $K$ at
$v$, as analog of $p$-adic situation, where  $p$ is a prime in
$\mathbb{Z}$ and $\mathbb{Q}_p$ the field of $p$-adic numbers.

For $v$ a prime of $\mathbb{F}_{q}[t]$, the $v$\emph{-adic zeta
function of Goss} is defined on the space
$S_{v}:=\mathbb{C}_{v}^{\ast}\times\underset{\overleftarrow{j}}{\lim
}\mathbb{Z}/(q^{\deg v}-1)p^{j}\mathbb{Z}$, where $\mathbb{C}_{v}$
is the completion of an algebraic closure of $K_v$ and the $\underset{\overleftarrow{j}}{\lim}\mathbb{Z}/(q^{\deg v}%
-1)p^{j}\mathbb{Z}$ is isomorphic to the product of $\mathbb{Z}_{p}$
with the cyclic group $\mathbb{Z}/(q^{\deg v}-1)\mathbb{Z}$. See
\cite[\S $8.3$]{bsffa}, \cite[\S $5.5$(b)]{ffa} for motivation and
details.

For $s=(x,y)\in S_{v}$, the Goss $v$-adic zeta
function is then defined as
\[
\zeta_{v}(s)=\sum_{d=0}^{\infty}x^{d}\sum_{_{\substack{a \in
\mathbb{F}_{q}[t],\text{ monic}
\\\deg(a)=d \\(v,a)=1}}}a^{y}.
\]
Note that $a^y$ here, by Fermat's little theorem, is a $p$-adic power of a
one-unit at $v$, and thus makes sense.
For $y\in \mathbb{Z}/(q^{\deg v}-1)\mathbb{Z} \times \mathbb{Z}_{p}$
we write
\[
S_{d,v}(y):=\sum_{\substack{a \in \mathbb{F}_{q}[t], \text{ monic}
\\\deg(a)=d
\\(v,a)=1}}a^{y},
\]
and let
\[
v_{d}(y):=\text{ val}_{v}(S_{d,v}(y)).
\]

We focus on $v$ of degree one, so that without loss of generality we assume
that $v=t$.

In this case, the Riemann hypothesis is the statement that for a
fixed $y$ all zeros in $x$ of  $\zeta_{v}(s)$ are in
$\mathbb{F}_{q}((t))=K_t$ (inside much bigger field $\mathbb{C}_t$),
and they are simple zeros, i.e., with multiplicity $1$.

We note that for $q=2$, and $v=t$,  Wan
\cite{rhcpzf} proved this and for general $q$, Goss \cite[Prop.
$9$]{rhcolf} got the same result for $y$ in $(q-1)S_v$. This last
result can also be derived from the results obtained by Thakur in
\cite[see corollary $8$]{vpsgvzf}. We drop these restrictions.

Let $q$ be a power of a prime $p$. Let $n$ be any positive integer.
If we can write $n$ as a sum of positive integers $n_0+ \cdots+
n_d$, so that when we consider $n_0,\ldots,n_d$ in their base-$p$
expansions, there is no carry over of coefficients when we perform
the sum in base-$p$, we express this by writing $n = n_0 \oplus
\cdots \oplus n_d$. Of all such decompositions of $n$, choose those
such that $n_d \leq \cdots\leq n_0$. Having done so, we can order
such decompositions of $n$ by saying $n_0 \oplus \cdots \oplus n_d
\leq m_0 \oplus \cdots \oplus m_d$ if for some $i,
n_d=m_d,\ldots,n_i=m_i$, but $n_{i-1}<m_{i-1}$, and will say that
the minimal such decomposition is ``given by the greedy algorithm.''

The next theorem is a slight generalization of Theorem $7$ (iii) in
\cite{vpsgvzf}, which is the special case $z \equiv m \mod (q-1)$.

\begin{teo}
\label{t 1}Let $m$ be a positive integer and $y=(z,m) \in \mathbb{Z}/(q^{\deg v}-1)\mathbb{Z} \times \mathbb{Z}_{p}$. Let $v$ be a prime of $\mathbb{F}%
_{q}[t]$ of degree one. Then either
$v_{d}(y)=\min(m_{1}+\cdots+dm_{d})$, where
$m=m_{0}\oplus\cdots\oplus m_{d}$, $(q-1)$ divides $m_{i}>0$ for
$0<i<d$, and $q-1$ divides $z-m_d$; or $v_{d}(k)$ is infinite, if
there is no such decomposition. When the decomposition exists, the
minimum is uniquely given by the greedy algorithm. (If the least
non-negative residue mod $q-1$ of $z$ is $r,m_{d}$ is the least
possible sum of $p$-powers chosen from the $p$-expansion of $m$
which is $r$ modulo $q-1$.)
\end{teo}

\proof

For $v=t$, we have
\begin{align*}
S_{d,v}(y)  &  =\sum_{\substack{f_{i}\in\mathbb{F}_{q}
\\f_{0}\neq0}} f_{0}^z \left(
1+\frac{f_1}{f_0}t+\cdots+\frac{f_{d-1}}{f_{0}}t^{d-1}+\frac{t^{d}}{f_{0}}\right)^{m}\\
&  =\sum_{_{\substack{f_{i}\in\mathbb{F}_{q} \\f_{0}\neq 0}}}f_{0}^z
\sum_{_{\substack{m=m_0+\cdots+m_d \\ m_i \geq 0,
i=0,\ldots,d}}}\binom{m}{m_{0},\ldots,
m_{d}}\left(\frac{f_1}{f_0}\right)^{m_1}\cdots
\left(\frac{f_{d-1}}{f_{0}}\right)^{m_{d-1}}\left(\frac{1}{f_0}\right)^{m_d} t^{m_{1}+\cdots+dm_{d}}\\
&  =\sum_{_{\substack{f_{i}\in\mathbb{F}_{q}
\\f_{0}\neq0}}}\sum_{_{\substack{m=m_0+\cdots+m_d \\ m_i \geq 0,
i=0,\ldots,d}}}\binom{m}{m_{0},\ldots, m_{d}}(f_{0})^{z-m_{d}}(f_{1})^{m_1}%
\cdots(f_{d-1})^{m_{d-1}}t^{m_{1}+\cdots+dm_{d}}.
\end{align*}
Now, by the theorem of Lucas $\binom{m}{m_{0},\ldots, m_{d}}$ is
non-zero if and only if the sum $m=m_{0}+\cdots+m_{d}$ is such that
there is no carry over of digits base $p$. Changing the order of the
summation and using the relation
\[
\sum_{f\in\mathbb{F}_{q}}f^{h}=\left\{
\begin{array}
[c]{ll}%
0 & \text{if }(q-1)\nmid h\text{ or }h=0\\
-1 & (q-1)\mid h,
\end{array}
\right.
\]
we see that
\[
S_{d,v}(y)=\sum_{_{\substack{f_{i}\in\mathbb{F}_{q} \\f_{0}\neq0}}}%
\sum_{\substack{m=m_{0}\oplus\cdots\oplus m_{d} \\(q-1)\mid m_{i}%
>0,i=1,\ldots,d-1 \\(q-1)\mid z-m_d \\ m_{0},m_d\geq0}}\binom{m}{m_{0},\ldots, m_{d}}%
(f_{0})^{z-m_{d}}(f_{1})^{m_1}\cdots(f_{d-1})^{m_{d-1}}t^{m_{1}+\cdots+dm_{d}}.
\]
By Sheats \cite{rhgzf}, $\min(m_{1}+\cdots+dm_{d})$, where $m=m_{0}%
\oplus\cdots\oplus m_{d}$, $(q-1)$ divides $m_{i}>0$ for $0<i<d$,
and $q-1$ divides $z-m_{d}$ is unique. From this, it follows that
$v_{d}(y)=\min(m_{1}+\cdots+dm_{d})$ is this unique
minimum.$\hfill\square$

\begin{com}
\label{c 1}We can go further in Theorem \emph{\ref{t 1}}. Not only
choose $m_{d}$ that way, but also choose $m_{d-1}$ to be the least
possible sum of $p$-powers chosen from the expansion of $m-m_{d}$
which is $0$ $\operatorname{mod}$ $q-1$, and so on.
\end{com}

In the case that $q=p$, a prime, we have the following special case
of Theorem \ref{t 1} (see (\cite[Thm. $9$]{vpsgvzf})).

\begin{teo}\label{t 2}
Let $q=p$ be a prime, $v$ a prime of $\mathbb{F}_{q}[t]$ of degree
one, $y=(z,m) \in \mathbb{Z}/(q^{\deg v}-1)\mathbb{Z} \times
\mathbb{Z}_{p}$ and $m>0$ an integer. Also, let $r$ be the least
non-negative residue of $z \mod (q-1) $. Write
$m=\sum_{i=1}^{\ell}p^{e_{i}}$, with $e_{i}$ monotonically
increasing and with not more than $p-1$ of the consecutive values
being the same (i.e., consider the base $p$-digit expansion
sequentially one digit at a time). Then $v_{d}(y)$ is infinite if
$\ell<(p-1)(d-1)+r$, and otherwise
\[
v_{d}(y)=d\sum_{i=1}^{r}p^{e_{i}}+\sum_{j=1}^{d-1}j\sum_{i=1}^{p-1}%
p^{e_{(d-1-j)(p-1)+r+i}}.
\]

\end{teo}

\proof

We have $v_{d}(y)=\min(m_{1}+\cdots+dm_{d})$, where
$m=m_{0}+\cdots+m_{d}$, when this sum has no carry over base $p$ and
$q-1$divides $m_{u}>0$ for $0<u<d$ and $q-1$ divides $z-m_{d}$. Note
that $p^{u}\equiv 1 \mod (q-1)$ when $q$ is prime. Hence $p-1$
powers together give divisibility by $q-1$. Hence for the minimum we
choose $r$ digits for $m_{d}$ and $m_{d-1},\cdots,m_{1}$ are
obtained by picking $p-1$ digits from the base $p$ expansion of
$m-m_{d}$ starting from the lowest digits (and dumping the rest of
the expansion, if any, into $m_{0}$).$\hfill\square$

With respect to the location of the zeros of the Goss $v$-adic zeta
function and its multiplicity, when the degree of $v$ is one, Wan
\cite{rhcpzf} showed that for $q = 2$, the zeros are simple with
corresponding $x$ in $\mathbb{F}_{q}((t))$. For general $q$, Goss
\cite[Prop. $9$]{rhcolf} got the same result for $y$ in $(q-1)S_v$.
This last result can also be derived from the results obtained by
Thakur in \cite[see corollary $8$]{vpsgvzf}. Here we prove the full
Riemann hypothesis for the $v$-adic zeta function for $v$ of degree
one, that is, all the zeros of $\zeta_v(x, y)$ are simple with  $x$
in $\mathbb{F}_{q}((t))$ inside $\mathbb{C}_t$.

For $k>0$, let $\ell(k)$ be the sum of the digits of the base $q$
expansion of $k$. Now, it is easy to show that, if $y=(z,m) \in
\mathbb{Z}/(q^{\deg v}-1)\mathbb{Z} \times \mathbb{Z}_{p}$, $m>0$ an
integer, and $d>\dfrac{\ell(m)}{q-1}+1$, then $S_{d,v}(y)=0$ (see
\cite[Thm. $5.1.2$]{ffa}). So, $\zeta_{v}(s)$ is a polynomial on
$x$. The Newton polygon of this polynomial determines the $v$-adic
valuation of its zeros. Namely, if the Newton polygon has a side of
slope $-m$ whose horizontal projection is of length $l$, then the
polynomial has precisely $l$ roots with valuation $m$. Furthermore,
if

\[\zeta_{m}(x)=\prod_{\stackrel{\zeta_{v}(\alpha)=0}{\text{val}_{v}(\alpha)=m}}(x-\alpha)\]
then $\zeta_{m}(x)\in \mathbb{F}_{q}((t))[x]$. Now, as in \cite{rh},
we have

\begin{lem}\label{l 1} Let $q=p^{n},n\geq1$. For $v=t$ prime of $\mathbb{F}_{q}[t]$ of degree one, and
for $y=(z,m) \in \mathbb{Z}/(q-1)\mathbb{Z} \times
\mathbb{Z}_{p}$ with $m\geq 0$ an integer , the zeros of
$\zeta_{v}(s)=\sum _{d=0}^{\infty}S_{d,v}(y)x^{d}$ are simple with $x$
 in $\mathbb{F}_{q}((t))$.
\end{lem}

\proof

Let us assume that Goss $t$-adic zeta function is a polynomial of
degree $d$,
\[\zeta_{v}(s)=1+S_{1,v}(y)x+ \cdots +S_{d,v}(y)x^{d},\]
and recall that $v_d(y):=\text{ val}_{v}(S_{d,v}(y))$. Then,
\[v_d(y)=m_1(d)+2m_2(d)+ \cdots +dm_d(d),\]
with obvious notation.

If $r$ is the least non-negative residue of $z$ modulo $q-1$, by
Theorem \ref{t 1} and Comment \ref{c 1}, $m_d (d)$ is the least
possible sum of $p$-powers chosen from the $p$-expansion of $m$
which is $r$ modulo $q-1$. Next, $m_{d-1}(d)$ is obtained as the
least possible sum of $p$-powers chosen from the expansion of
$m-m_d(d)$ which is $0$ mod $q-1$, and so on.

For $q=p$, a prime number, the greedy solution is achieved in a
simpler way, by the previous Theorem \ref{t 2}:
$m_{d}(d)=\sum_{s=1}^{r}p^{e_{s}},m_{d-1}(d)=\sum_{s=1}^{p-1}p^{e_{r+s}}$,
and so on.

Therefore, for $q$ a prime or power of a prime, the valuations
$v_d(y)$ of the coefficients of Goss $v$-adic zeta functions are:

\begin{eqnarray*}
  v_d(y)&=& m_1(d)+2m_2(d)+ \cdots +dm_d(d) \\
  v_{d-1}(y)&=& m_1(d-1)+2m_2(d-1)+ \cdots +(d-1)m_{d-1}(d-1) \\
  \vdots&=& \vdots \\
  v_1(y)&=& m_1(1) \\
  v_0(y) &=& 0.
\end{eqnarray*}

Now, notice that by construction we have
\begin{equation*}
m_{d-j}(d)=m_{d-j-1}(d-1)\text{ for }j=0,\cdots,d-2.
\end{equation*}
Using this, straight manipulation allows us to calculate the slope
$\lambda(d)$ of the Newton polygon from
$(d-1,v_{d-1}(y))$ to $(d,v_{d}(y)),$%

\[
\lambda(d)=v_{d}(y)-v_{d-1}(y)=m_{1}(d)+m_{1}(d-1)+m_{2}(d-1)+\cdots
+m_{d-1}(d-1).
\]
Finally, $\lambda(d)$ is a strictly increasing function of $d$ since
\[
\lambda(d)-\lambda(d-1)=m_{1}(d)>0
\]
for $d>1$. For $d=1, \lambda(d)-\lambda(d-1)=m_{1}(1)$, if we put
$\lambda(0)=0$, and in this case $m_{1}(1)$ could be zero if $z$ is
even, i.e., when $z$ is congruent to $0$ mod $q-1$, but that just
states that first slope is zero. This implies that each slope of the
Newton polygon has horizontal projection one and so the zeros of
$\zeta_{v}(s)$ are simple with $x$ in $\mathbb{F}_{q}((t))$.
 $\hfill\square$

\begin{teo}[Riemann Hypothesis] Let $q$ be any prime power.
For the prime $v=t$ of $\mathbb{F}_{q}[t]$, and for $y=(z,M) \in
\mathbb{Z}/(q-1)\mathbb{Z} \times \mathbb{Z}_{p}$, the zeros of
$\zeta_{v}(x,y)$ are simple and have the  $x$  in
$K_v=\mathbb{F}_{q}((t))$ inside $\mathbb{C}_v$.
\end{teo}

\proof By the Lemma 1, we can assume that  $M \in \mathbb{Z}_{p}$
 is not a non-negative integer. Thus we can
choose sufficiently many  appropriate $p$-adic digits to form a positive integer $m> 0$,
which has a decomposition as in Theorem \ref{t 1} with positive
$m_i$'s for all $i=1, \cdots, d-1$, so that the valuation stabilizes
(exactly as in  \cite{rh}, \cite[Sec. 5.8]{ffa}), as it is independent of
$m_0$. We can form pairs $ y_m=(z, m),m>0$ that converge to $(z,
M)$. Hence we can calculate the valuation as
\begin{eqnarray*}
v_d(y) &=&\text{ val}_{v}\left(\sum_{_{\substack{a \in \mathbb{F}_q[t],\text{monic}\\\deg(a)=d \\(v,a)=1}}}a^{y}\right) \\
       &=&\text{ val}_{v}\left(\sum_{_{\substack{a \in \mathbb{F}_q[t], \text{monic}\\\deg(a)=d \\(v,a)=1}}}a^{y_m}\right) \\
       &=&v_d(y_m).
\end{eqnarray*}
Then  the slope inequalities  of the Lemma \ref{l 1} imply the required result in
general. $\hfill\square$

\textbf{Acknowledgement}

We are grateful to Dinesh Thakur for having suggested us this
project. We also thank the referee for the thoughtful comments and
suggestions.


\begin{thebibliography}{9}

\bibitem {rh}J. Diaz-Vargas, Riemann hypothesis for $\mathbb{F}_{q}[t]$, J.
Number Theory $\mathbf{59}$ $(1996),313-318$.

\bibitem {bsffa}D. Goss, \emph{Basic Structures of Function Field Arithmetic},
Springer-Verlag $1998$.

\bibitem{rhcolf}D. Goss, A Riemann hypothesis for characteristic $p$
 $L$-functions, J. Number Theory, $\mathbf{82}$ $(2000),
299-322$.

\bibitem {rhgzf}J. Sheats, The Riemann hypothesis for the Goss zeta function
for $\mathbb{F}_{q}[t]$, J. Number Theory $\mathbf{71}$ $(1998),121-157$.

\bibitem {ffa}D. Thakur, \emph{Function Field Arithmetic}, World Scientific,
$2004$.

\bibitem {vpsgvzf}D. Thakur, Valuations of $v$-adic power sums and zero
distribution for the Goss $v$-adic zeta function for $\mathbb{F}_{q}[t]$,
\emph{J. Integer Seq.}, $\mathbf{16}$ $(2013)$, pp. 1-18.  Article $13.2.13$.

\bibitem {rhcpzf}D. Wan, On the Riemann hypothesis for the characteristic
$p$ zeta function, J. Number Theory, $\mathbf{58}$ $(1996),
196-212$.

\end{thebibliography}
\end{document}